\input amstex
\documentstyle{amsppt}
\NoRunningHeads
\NoBlackBoxes
\document

\def\Q{{\Bbb Q}}

\def\ell{{\text{ell}}}

\def\h1{\hat{\bold 1}}

\def\a{\frak a}

\def\Ua{U_q(\tilde\g)}
\def\U2{{\Ua}_2}
\def\g{\frak g}
\def\n{\frak n}

\def\Z{\Bbb Z}
\def\C{\Bbb C}

\def\<{\langle}
\def\>{\rangle}
\def\o{\otimes}

\def\b{\frak b}
\def\h{{\frak h}}
\topmatter
\title Quantization of Lie bialgebras, part VI: 
quantization of generalized Kac-Moody algebras
\endtitle
\author {\bf Pavel Etingof
\footnote{Department of Mathematics, 2-165, MIT,
Cambridge, MA 02139, USA, etingof\@math.mit.edu} and David Kazhdan
\footnote{Einstein Institute of Mathematics
Edmond J. Safra Campus, Givat Ram
The Hebrew University of Jerusalem
Jerusalem, 91904, Israel, kazhdan\@math.huji.ac.il}}
\endauthor
\endtopmatter

\centerline{\bf To Bert Kostant with admiration}

\heading 1. Introduction.\endheading

\vskip .1in

This paper is a continuation of the series \cite{EK1-5}.
We show that the image of a Kac-Moody Lie bialgebra 
with the standard quasitriangular structure 
under the quantization functor defined in \cite{EK1,EK2}
is isomorphic to the Drinfeld-Jimbo quantization 
of this Lie bialgebra, with the standard quasitriangular
structure. This implies that 
when the quantization parameter is formal, then the category 
$\Cal O$ for the quantized Kac-Moody algebra is equivalent, 
as a braided tensor category, to the category $\Cal O$ 
over the corresponding classical Kac-Moody algebra,
with the tensor category structure defined by a Drinfeld
associator. This equivalence is a generalization of the functor
constructed in \cite{KL}. 

In particular, we answer positively questions 8.1, 8.2
from \cite{Dr1}: we show that the characters of irreducible
highest weight modules for a quantized Kac-Moody algebra 
$\g$ are the same as in the classical case, and that the quantum
deformation of the appropriate completion $\hat U(\g)$ of
$U(\g)$ is trivial as a deformation of algebras. 

Moreover, our results are valid for the Lie algebra 
$\g(A)$ corresponding to any symmetrizable matrix $A$
(not necessarily with integer entries). 
This answers question 8.3 in \cite{Dr1}, of Drinfeld and 
Gelfand (how to define a flat deformation 
$U_\hbar(\g(A))$ of the Hopf algebra $U(\g(A))$ 
for any symmetrizable matrix $A$).  

We also prove the Drinfeld-Kohno theorem 
for the algebra $\g(A)$ (it was previously proved by Varchenko 
\cite{V} using integral formulas for solutions of the KZ equations).

{\bf Remark.} One of the important facts used in this paper 
is that the quantization functors from
\cite{EK1,EK2} commute with duals and doubles.
However, the original version of this paper, which 
appeared in 2000, unfortunately did not contain a convincing proof of 
this statement. Namely, it referred 
to \cite{EK1,EK2,EK3}, where this fact was proved for finite
dimensional Lie bialgebras, and claimed that in general the proof
was similar. But it turns out that actually, it is not easy
to extend the argument of \cite{EK1,EK2,EK3} to the general case,
since it uses cyclic expressions. A general proof of this result 
was obtained by Enriquez and 
Geer in 2007, see \cite{EG}. This revised version of our paper 
takes this fact into account, and replaces insufficient references
to \cite{EK1,EK2,EK3} by the reference to \cite{EG}. 

{\bf Acnowledgments.} The authors were supported by the NSF
grant DMS-9700477. 

\heading 2. Generalized Kac-Moody algebras\endheading 

Throughout the paper, $k$ denotes a field of characteristic
zero. All vector spaces in this paper will be over $k$. 

{\bf 2.1.} We recall definitions from \cite{K}. 
Let $A=(a_{ij})$ be an $n$-by-$n$ matrix with entries in $k$, 
and $(\h,\Pi,\Pi^\vee)$ be a realization of $A$. 
This means that 
$\h$ is a vector space of dimension $2n-\text{rank}(A)$,
$\Pi=\{\alpha_1,...,\alpha_n\}\subset \h^*$, 
$\Pi^\vee=\{h_1,...,h_n\}\subset \h$ are
linearly independent, and $\alpha_i(h_j)=a_{ji}$.

\proclaim{Definition} The Lie algebra $\tilde \g(A)$ is generated 
by $\h,e_1,...,e_n,f_1,...,f_n$ with defining relations 
$$
[h,h']=0,\ h,h'\in\h;\ [h,e_i]=\alpha_i(h)e_i; \ 
[h,f_i]=-\alpha_i(h)f_i;\ [e_i,f_j]=\delta_{ij}h_i.
$$
\endproclaim

We will denote $\tilde \g(A)$ simply by $\tilde \g$, 
assuming that $A$ has been
fixed.

Let $I$ be the sum of all two-sided ideals in $\tilde\g(A)$ which 
have zero intersection with $\h\subset \tilde \g(A)$. 
Let $\g(A):=\tilde\g(A)/I$. 
The algebra $\g(A)$ is called a generalized Kac-Moody algebra. 
We will denote $\g(A)$ by $\g$, assuming that $A$ has been
fixed. The Lie algebra $\g$ is graded by principal gradation
($\text{deg}(e_i)=1,\text{deg}(f_i)=-1,\text{deg}(\h)=0$), and 
the homogeneous components are finite dimensional. 

In the following we will assume that the matrix $A$ is 
symmetrizable, i.e. there exists a collection of nonzero numbers 
$d_i$, $i=1,...,n$, 
such that $d_ia_{ij}=d_ja_{ji}$. We will choose such a collection of
numbers. 
Let us choose a nondegenerate 
bilinear symmetric form on $\h$ such that 
$(h,h_i)=d_i^{-1}\alpha_i(h)$.
It is easy to see that such a form always exists.  
It is known \cite{K} that there exists a unique extension of the
form $(,)$ to an invariant symmetric bilinear form
$(,)$ on 
$\tilde \g$. (For this extension, one has 
$(e_i,f_j)=\delta_{ij}d_i^{-1}$). The kernel of this form is 
$I$, and thus the form descends to a nondegenerate form 
on $\g$. 

{\bf Remark.}
One can show that forms on $\g$ coming from different
forms on $\h$ are equivalent under automorphisms of $\g$. 

{\bf 2.2.} Let $\n_+,\n_-,\b_+,\b_-$ be the nilpotent and the
Borel subalgebras of $\g$ ($\n_+,\n_-$ are generated by 
$e_i$ and by $f_i$, respectively, and $\b_\pm:=\n_\pm\oplus \h$).
Let us regard $\b_+$ and $\b_-$ as Lie subalgebras 
of $\g\oplus \h$ using the embeddings 
$\eta_\pm: \b_\pm \to \g\oplus \h$ given by
$$
\eta_\pm(x)=x\oplus (\pm \bar x),
$$
where $\bar x$ is the image of $x$ in $\h$. 

Define the inner product on $\g\oplus\h$ by $(,)_{\g\oplus \h}=
(,)_\g-(,)_\h$. 
The following proposition is well known and straightforward to
check. 

\proclaim{Proposition 2.1} The triple $(\g\oplus\h, \b_+,\b_-)$
with the inner product $(,)_{\g\oplus \h}$ 
and embeddings $\eta_\pm$ is a (graded) Manin
triple. 
\endproclaim

The proposition implies that $\g\oplus \h$, $\b_+$, $\b_-$ are naturally Lie
bialgebras, with $\b_+^*=\b_-^{op}$, where ${}^*$ denotes the
restricted dual space, and $ ^{op}$ denotes 
the opposite cocommutator. The cocommutator $\delta$ on these algebras
is easily computed:
$$
\delta(h)=0, h\in \h\subset \b_\pm;\ 
\delta(e_i)=\frac{1}{2}d_ie_i \wedge h_i;\ \delta(f_i)=
\frac{1}{2}d_if_i\wedge h_i.\tag 2.1
$$

The Lie subalgebra $\{(0,h)|h\in \h\}$ is thus an ideal and
coideal in $\g\oplus \h$, and so the quotient $\g=(\g\oplus\h)/\h$
is also a Lie bialgebra with Lie subbialgebras $\b_+,\b_-$, 
and the same cocommutator formulas. 

In fact, the same formulas define a Lie bialgebra structure
on $\tilde \g$ and its Borel subalgebras 
$\tilde \b_\pm$ (generated by $\h,e_i$ and $\h,f_i$,
respectively). The projections 
$\tilde \g\to \g$, $\tilde \b_\pm\to \b_\pm$ are thus 
Lie bialgebra homomorphisms. 

{\bf Remark.} The factors $\frac{1}{2}$ in (2.1) appear because 
for $a\in \h\subset \b_+$ and $b\in \h\subset \b_-$, 
one has $(a,b)_{\g\oplus \h}=(a+\bar a,b-\bar b)
=2(a,b)$.

\heading 3. Quantization of generalized Kac-Moody
algebras\endheading 

{\bf 3.1.} Let $\hbar$ be a formal parameter, and
$q=e^{\hbar/2}$.  
Let $\Phi=1+\frac{\hbar^2}{24}[\Omega_{12},\Omega_{23}]+...$ 
be a universal Lie associator (see \cite{Dr4}), and $\a\to U_\hbar(\a)$ 
be the functor of quantization of Lie bialgebras associated 
with $\Phi$ (see \cite{EK1,EK2}). In this section 
we will describe explicitly $U_\hbar(\a)$, when $\a$ is one 
of the Lie bialgebras of the previous section. 

\proclaim{Proposition 3.1} The QUE algebra $U_\hbar(\tilde \b_+)$ 
is isomorphic to the QUE algebra $\tilde {\Cal U}_+$ generated
(topologically) by $\h$ and elements
$E_i$, $i=1,...,n$, with the relations 
$$
[h,h']=0;\  [h,E_i]=\alpha_i(h)E_i, h,h'\in \h,
$$
with coproduct 
$$
\Delta(h)=h\otimes 1+1\otimes h;\ \Delta(E_i)=E_i\o q^{\gamma_i}+1\otimes E_i,
$$
for suitable elements $\gamma_i\in \h[[\hbar]]$.  
\endproclaim

\demo{Proof} Since $U_\hbar$ is a functor, the embedding of Lie
bialgebras $\h\to \tilde \b_+$ defines an embedding of 
QUE algebras $U_\hbar(\h)=U(\h)[[\hbar]]\to
U_\hbar(\tilde\b_+)$. 

Also, $\b_+$ has a $\Z_+^n$-grading given by 
$\text{deg}_i(h)=0$, $\text{deg}_i(e_j)=\delta_{ij}$, so 
by functoriality the quantized algebra $U_\hbar(\tilde\b_+)$ 
has a grading by $\Z_+^n$ as well (as this grading is simply 
an action of $\Bbb G_m^n$). 

As a result, we get $U_\hbar(\tilde\b_+)=
\oplus_{\bold m\in \Bbb Z_+^n}U_\hbar(\tilde\b_+)[\bold m]$, 
where $U_\hbar(\tilde\b_+)[\bold m]$
are free modules over $U_\hbar(\h)$ of finite rank
(in fact, the same rank as before deformation). 
In particular, if $\bold m=1_j$, where $1_j(i)=\delta_{ij}$, 
then $U_\hbar(\tilde\b_+)[\bold m]$ has rank $1$. 

Let us choose an element $E_j'$ in $U_\hbar(\tilde\b_+)[1_j]$
which equals $e_j$ modulo $\hbar$. 

For homogeneity reasons we have 
$$
\Delta(E_j')=(E_j'\otimes 1)\Psi_1
+(1\otimes E_j')\Psi_2,
$$
where $\Psi_i\in 1+\hbar U(\h\oplus \h)[[\hbar]]$. 

We have $(\Delta\o 1)(\Delta(E_i'))=
(1\o \Delta)(\Delta(E_i'))$. This implies the following equations 
on $\Psi_1,\Psi_2$:
$$
(\Psi_1\otimes 1)(\Delta\otimes
1)(\Psi_1)=(1\otimes\Delta)(\Psi_1),\tag 3.1
$$
$$
(\Psi_2\otimes 1)(\Delta\o 1)(\Psi_1)=
(1\o \Psi_1)(1\o \Delta)(\Psi_2),\tag 3.2
$$
$$
(1\o \Psi_2)(1\o \Delta)(\Psi_2)=(\Delta\o 1)(\Psi_2).\tag 3.3
$$

Let us regard $\Psi_i$ as functions of two variables 
$x,y\in \h^*$, and let $\psi_i=\log \Psi_i$. Then (3.1)-(3.3) 
can be written in the form
$$
\psi_1(x,y)+\psi_1(x+y,z)=\psi_1(x,y+z),\tag 3.4
$$
$$
\psi_2(y,z)+\psi_2(x,y+z)=\psi_2(x+y,z),\tag 3.5
$$
$$
\psi_1(y,z)+\psi_2(x,y+z)=\psi_1(x+y,z)+\psi_2(x,y).\tag 3.6
$$
Let us set $z=-x-y$ in (3.4). We get
$$
\psi_1(x,y)=\phi_1(x)-\phi_1(x+y),\tag 3.7
$$
where $\phi_1(x)=\psi_1(x,-x)$. 
Similarly, from equation (3.6), putting $x=-y-z$, we get 
$$
\psi_2(y,z)=\phi_2(z)-\phi_2(y+z),\tag 3.8
$$
where $\phi_2(z)=\psi_2(-z,z)$. 

It is easy to check that 
after these substitutions, equation (3.5) becomes
$$
\phi_1(y)-\phi_1(x+y)-\phi_1(y+z)+\phi_1(x+y+z)=
\phi_2(y)-\phi_2(x+y)-\phi_2(y+z)+\phi_2(x+y+z).\tag 3.9
$$
Let $\gamma(x)=\frac{1}{\log q}(\phi_1(x)-\phi_2(x))$. We have 
$$
\gamma(y)-\gamma(x+y)-\gamma(y+z)+\gamma(x+y+z)=0.\tag 3.10
$$
In particular, $d^2\gamma=0$ and hence $\gamma$ is an affine linear
function (i.e. $\gamma\in  (\h\oplus k)[[\hbar]]$).
We will denote $\gamma$ by $\gamma_i$ to remember 
its dependence on $i$. 

Define $E_i=E_i'e^{-\phi_2(x)}$. Then it is easy to see 
from the above that 
$$
\Delta(E_i)=E_i\otimes q^{\gamma_i}+1\otimes E_i.\tag 3.11
$$
 From this we see using the counit axiom that 
the constant terms of $\gamma_i$ are zero, 
i.e. $\gamma_i$ are elements of $\h[[\hbar]]$.  

It is also clear that $\h$ and $E_i$ topologically generate
$U_\hbar(\tilde\b_+)$ and that the only relations are the ones 
given in the theorem (this follows from the fact that $\tilde
\n_+$ is a free Lie algebra). The proposition is proved. 
$\square$
\enddemo

{\bf 3.2.} Let us now compute the elements $\gamma_i$. 

\proclaim{Proposition 3.2} One has $\gamma_i=d_ih_i$. 
\endproclaim

\demo{Proof} We have a surjective map of Lie bialgebras
$\tilde \b_+\to \b_+$. By functoriality of quantization, 
this defines a surjective 
homomorphism of QUE algebras $U_\hbar(\tilde \b_+)\to U_\hbar(\b_+)$
(which preserves the grading). Therefore, 
$U_\hbar(\b_+)$ is also generated 
by $\h,E_i$ satisfying the relations 
$$
[h,E_i]=\alpha_i(h)E_i
$$
(and maybe some additional relations), and 
the coproduct is defined by 
$$
\Delta(h)=0,\ \Delta(E_i)=E_i\otimes q^{\gamma_i}+1\otimes E_i.
$$

It follows from the definition of the Lie bialgebra $\b_+$ 
that it is self-dual in the graded sense: $\b_+\cong \b_+^*$.
Thus, $U_\hbar(\b_+)\cong U_\hbar(\b_+^*)$. 
By the result of \cite{EG}, we have 
$U_\hbar(\b_+^*)\cong U_\hbar(\b_+^{op})^{*op}$. 
On the other hand, for any Lie bialgebra $\a$, 
$U_\hbar(\a^{op})\cong U_{-\hbar}(\a)$ (since the universal quantization
formulas of \cite{EK2} are written in terms of $[,]$ and $\hbar\delta$).
Therefore, $U_\hbar(\b_+^{op})^{*op}\cong 
U_{-\hbar}(\b_+)^{*op}$. Thus, we have an (graded) isomorphism 
of QUE algebras $U_\hbar(\b_+)\to U_{-\hbar}(\b_+)^{*op}$,
which in degree zero comes from the identification 
$\h\to \h^*$ using the form $2(,)$ on $\h$
. This isomorphism can be
understood as a bilinear form 
$B:U_\hbar(\b_+)\otimes U_{-\hbar}(\b_+)\to k((\hbar))$ 
satisfying the conditions 
$$
B(xy,z)=B(y\otimes x,\Delta(z)), \ 
B(z,xy)=B(\Delta(z),x\otimes y),
$$
such that 
$$
B(q^{a},q^{b})=q^{-(a,b)}, a,b\in \h.
$$
Let $B_i=B(E_i,E_i)$; clearly, this is nonzero.
Using the properties of $B$, we have 
$$
B(E_i,E_iq^{ a})=q^{-(a,\gamma_i)}B_i;\
B(E_i,q^{a}E_i)=B_i.
$$
But $q^aE_iq^{-a}=q^{\alpha_i(a)}E_i$, so we get
$(a,\gamma_i)=\alpha_i(a)$, which yields $\gamma_i=d_ih_i$, as
desired. 
$\square$\enddemo

Thus we have proved 

\proclaim{Theorem 3.3}
The QUE algebra $U_\hbar(\tilde \b_+)$ 
is isomorphic to the QUE algebra $\tilde {\Cal U}_+$ generated
(topologically) by $\h$ and elements
$E_i$, $i=1,...,n$, with the relations 
$$
[h,h']=0,\  [h,E_i]=\alpha_i(h)E_i, h,h'\in \h,
$$
with coproduct 
$$
\Delta(h)=0, \Delta(E_i)=E_i\o q^{d_ih_i}+1\otimes E_i.
$$
\endproclaim

{\bf 3.3.} Now let us describe explicitly the QUE algebra
$U_\hbar(\b_+)$. 

\proclaim{Theorem 3.4} There exists a unique 
symmetric bilinear form $B$ on $U_\hbar(\tilde \b_+)$ 
with values in $k((\hbar))$ which satisfies the properties 
$$
B(xy,z)=B(x\otimes y,\Delta(z)), \ 
B(z,xy)=B(\Delta(z),x\otimes y),
$$
$$
B(q^{ a},q^{b})=q^{-(a,b)}, a,b\in \h,
$$
$$
B(E_i,E_j)=\frac{\delta_{ij}}{q-q^{-1}}.
$$
The QUE algebra $U_\hbar(\b_+)$ is 
isomorphic to the quotient ${\Cal U}_+$ 
of $\tilde {\Cal U}_+$ (as in Theorem 3.3) by 
the Hopf ideal $Ker(B)$. 
\endproclaim

\demo{Proof} The existence and uniqueness of $B$ easily follows
from the freeness of the algebra generated by $E_i$. 
Also, the uniqueness of $B$ implies that $B$ is symmetric. 

Moreover, we claim that 
a form with the same properties exists and is unique on 
the quotient algebra $U_\hbar(\b_+)$. 
Indeed, we have seen that $U_\hbar(\b_+)\cong U_{-\hbar}(\b_+)^{*op}$. 
On the other hand,  by using the conjugation 
by $q^{-\sum x_i^2/2}$, where $x_i$ is an
orthonormal basis of $\h$, we find that $U_{-\hbar}(\b_+)^{op}\cong
U_\hbar(\b_+)$. This yields a symmetric
isomorphism 
$U_\hbar(\b_+)\to U_\hbar(\b_+)^*$, which gives a desired
form. (Here we use that for any Hopf algebra $H$, 
the Hopf algebra $H^{*op}$ is isomorphic to
$H^{op*}$ via the antipode.) 
 
Thus, the form on the big algebra is pulled back from the form on
the small algebra, and hence the kernel of the projection 
from the big algebra to the small one is contained in $Ker(B)$. 

Finally, it is clear that the form $B$ on $U_\hbar(\b_+)$ is
nondegenerate
(as it corresponds to an isomorphism). 
So $U_\hbar(\b_+)=\tilde {\Cal U}_+/Ker(B)={\Cal U}_+$, as desired. 
$\square$
\enddemo

\proclaim{Corollary 3.5} 
${\Cal U}_+$ is a flat defomation of $U(\b_+)$. 
\endproclaim

\demo{Proof} Clear, as $U_\hbar(\b_+)$ is flat by definition.
$\square$\enddemo

\proclaim{Corollary 3.6} Suppose that $A$ is a generalized 
Cartan matrix (i.e. $a_{ii}=2$, and 
$a_{ij}$ are nonpositive integers for $i\ne j$). 
In this case, the two-sided ideal $Ker(B)$ is generated by the 
quantum Serre relations 
$$
\sum_{m=0}^{1-a_{ij}}\frac{(-1)^m}
{[m]_{q_i}![1-a_{ij}-m]_{q_i}!}E_i^{1-a_{ij}-m}E_jE_i^m=0,
$$
where $q_i=q^{d_i}$.
\endproclaim

\demo{Proof} It is known (\cite{L}, Section 1) that the ideal generated by 
the quantum Serre relations is contained in 
$Ker(B)$. Besides, we know \cite{K} that $\b_+$
is the quotient of $\tilde \b_+$ by the classical limits 
of the Serre relations. This fact and Corollary 3.5 imply the result. 
$\square$\enddemo

\proclaim{Theorem 3.7} The QUE algebra $U_{\hbar}(\g)$ is isomorphic 
to the quotient ${\Cal U}$ of the (restricted) quantum double
$D({\Cal U}_+)$ by the ideal generated by the identification of
$\h\subset {\Cal U}_+$ and $\h^*\subset {\Cal U}_+^*$. In particular, 
if $A$ is a generalized Cartan matrix then 
$U_\hbar(\g)$ is isomorphic to the Drinfeld-Jimbo quantum group 
associated to the Kac-Moody algebra $\g$ (see 
\cite{Dr2}, Example 6.2, and \cite{J}). 
\endproclaim 

{\bf Remark.} The word ``restricted'' means that 
as a $k[[\hbar]]$-module, \linebreak $D(\Cal U_+)=\Cal U_+\otimes 
\Cal U_+^*$, where $\Cal U_+^*$ is the restricted (by the
grading) dual space to $\Cal U_+$. 

\vskip .05in

Thus, Theorem 3.7 constructs a flat deformation of $U(\g)$, and
thus answers Question 8.3 from \cite{Dr1}.

\demo{Proof} This follows from the previous results and the fact that 
quantization commutes with taking the double, see \cite{EG}.
$\square$ \enddemo

\vskip .01in

Now define $\tilde \g'$ to be the (restricted) Drinfeld
double of $\tilde\b_+$, as a Lie bialgebra. 

\vskip .01in

{\bf Remark.} We note that while for generic $A$ we have
$\tilde\g\cong\tilde\g'$ as graded Lie algebras, 
for special values of $A$ this is not
the case, and in particular the Lie
algebra
$\tilde\g'$ is not generated by elements of degree $1$ and $-1$. 

\proclaim{Theorem 3.8} 
The QUE algebra $U_{\hbar}(\tilde\g')$ is isomorphic 
to the quotient $\tilde {\Cal U'}$ of the (restricted) quantum double
$D(\tilde {\Cal U}_+)$ by the ideal generated by the identification of
$\h\subset \tilde{\Cal U}_+$ and $\h^*\subset \tilde 
{\Cal U}_+^*$. 
\endproclaim

\demo{Proof} The proof of this theorem is the same as 
the proof of Theorem 3.7. $\square$\enddemo

\heading 4. Category $\Cal O$\endheading

{\bf 4.1.} Let $\g_+$ be a 
Lie bialgebra, and $U_\hbar(\g_+)$ be its quantization as in
\cite{EK2}. Then one can define the standard notion of a Drinfeld-Yetter
module, or a dimodule, over $\g_+$ and $U_\hbar(\g_+)$
(see e.g. \cite{EK2}). Let $\Cal M$ be the category of deformation
dimodules over $\g_+$, i.e. of 
$\g_+$-dimodules realized on a topologically free $k[[\hbar]]$-module. 
Let $\Cal M_\hbar$ be the category of dimodules over
$U_\hbar(\g_+)$. Recall that both categories are braided 
tensor categories: the category $\Cal M$ has the braided tensor
structure defined by the associator $\Phi$, while 
the category $\Cal M_\hbar$ has the braided tensor structure 
obtained from the ``universal R-matrix'' (see \cite{EK2}).

\proclaim{Theorem 4.1}
There exists an equivalence of braided 
tensor categories $\Cal M\to \Cal M_\hbar$, 
which is the identity functor at the level of
$k[[\hbar]]$-modules
(i.e., there exists a consistent system of isomorphisms
of $k[[\hbar]]$-modules, $V\to F(V)$). 
\endproclaim

\demo{Proof} We will use the notation of \cite{EK1,EK2}.
Recall that in \cite{EK1,EK2}, we defined the functor 
$F$ from the category of deformation $\g_+$-dimodules
to the category of $k[[\hbar]]$-modules by 
$F(V):=\text{Hom}(M_-,M_+^*\hat\otimes V)$, and equipped it with a
tensor structure (here $\hat\otimes$ is the completed tensor
product with respect to the weak topology in $M_+^*$). To turn $F$ into a functor we are looking for,
we need to introduce on $F(V)$, for all $V$, an action and a
coaction of $U_\hbar(\g_+)$. 

Recall that $U_\hbar(\g_+)$ is defined in \cite{EK1,EK2} 
to be the space $F(M_-)$.    
The action of $U_\hbar(\g_+)$ on $F(V)$ is explicitly defined 
in \cite{EK1}, section 9. Namely, if $v\in F(V)$ and 
$a\in U_\hbar(\g_+)=F(M_-)$, one defines $av\in F(V)$ to be 
$(i_+^*\o 1\o 1)\circ (1\otimes v)\circ a$. 

Now define the coaction of $U_\hbar(\g_+)$ on $F(V)$. 
Since $F$ is a tensor functor, the braiding map 
for $\g_+$-dimodules composed with the permutation of components 
defines a map $R: F(M_-)\o F(V)\to F(M_-)\o
F(V)$ (the universal R-matrix). The coaction of $U_\hbar(\g_+)$ is defined 
by the map $v\to R(1\otimes v)$, where $1$ is the unit of
$U_\hbar(\g_+)$. One can check that these action and coaction 
are compatible, so they define a structure of a 
$U_\hbar(\g_+)$-dimodule on $F(V)$. Thus $F$ becomes a functor 
from $\g_+$-dimodules to $U_\hbar(\g_+)$-dimodules. 
It is straightforward to check that this functor 
equipped with the tensor structure of \cite{EK1} is a braided 
tensor functor between these categories. 

It remains to show that $F$ is an equivalence of categories. 
To do this, it is sufficient to construct the inverse functor. 
This is done using twisting the tensor category of dimodules 
of $U_\hbar(\g_+)$ by a family $a(t)$ of elements of the
Grothendieck-Teichm\''uller group, as explained in Section 2 of
\cite{EK2}. The theorem is proved.   
$\square$\enddemo

{\bf Remark.} We use this opportunity to correct the formulation 
of Theorem 6.2 in \cite{EK1}, whose original formulation is 
not quite correct. Instead of the category $\Cal
M_{\a}$ of $\a$-modules, considered in this theorem, 
one should consider the category $\tilde \Cal M_\a$ of deformation 
$\a$-modules. The functor $F$ in the theorem (from 
$\Cal M_\a$ to the category $\Cal R$ of representations of $U_\hbar(\a)$)
naturally 
extends to $\tilde \Cal M_\a$. The correct formulation 
of Theorem 6.2 says that $F$ is an equivalence of $\tilde \Cal M_\a$
onto $\Cal R$
(the proof of this is obvious from the results of \cite{EK1}). 
In this form, Theorem 6.2 of \cite{EK1} (for $\a$ being 
the double of a finite dimensional Lie bialgebra) is a special case of
Theorem 4.1 above. 
 
{\bf 4.2.} Let us return to the setting of generalized Kac-Moody
algebras. 
Recall that the category $\Cal O$ for $\g$ is defined 
to be the category of $\h$-diagonalizable 
$\g$-representations, whose weights belong to 
a union of finitely many cones $\lambda-\sum_i\Bbb Z_+\alpha_i$, 
$\lambda\in \h^*$, 
and the weight subspaces are finite dimensional. 
Define also the category $\Cal O[[\hbar]]$ of
deformation representations of $\g$, i.e. representations of 
$\g$ on topologically free $k[[\hbar]]$-modules with the above
properties (with $\lambda\in \h^*[[\hbar]]$).  

In a similar way one defines the category $\Cal O_\hbar$ 
for the algebra ${\Cal U}$: it is the category 
of ${\Cal U}$-modules which are topologically free over $k[[\hbar]]$ 
and satisfy the same conditions as in the classical case. 

Let $\Omega\in \g\hat\otimes\g$ 
(where $\hat\otimes$ is the tensor product completed with respect
to the grading)
be the inverse element to the bilinear form 
$(,)$ on $\g$. It defines an operator in any tensor product
$V\otimes W$ of
modules from category $\Cal O[[\hbar]]$. 
Following Drinfeld, we put on
$\Cal O[[\hbar]]$ a structure of a braided tensor category using the
associator $\Phi$, with braiding $q^{\Omega}$
(see \cite{EK1}). The category $\Cal O_\hbar$ is also a braided
tensor category, with braiding defined by the universal
R-matrix $R\in U_\hbar(\b_+)\hat\otimes U_\hbar(\b_-)$ 
coming from the isomorphism $U_\hbar(\b_+)\to
U_\hbar(\b_-)^{*op}$. 

Recall that a highest weight module over $\g$ or $U_\hbar(\g)$ 
is a module generated by a highest weight vector, and that for each
highest weight $\lambda$, we have the Verma module $M(\lambda)$
and the irreducible module $L(\lambda)$, and any highest weight
module $N$ with highest weight $\lambda$ can be included in 
a diagram $M(\lambda)\to N\to L(\lambda)$, where both maps are
surjective, and defined uniquely up to scaling.  

\proclaim{Theorem 4.2} There exists an equivalence of braided 
tensor categories $F: \Cal O[[\hbar]]\to \Cal O_\hbar$, 
which is isomorphic to the identity functor at the level of
$\h$-graded $k[[\hbar]]$-modules. This equivalence maps 
the Verma (resp. irreducible) module with highest weight
$\lambda$ to the 
the Verma (resp. irreducible) module with highest weight
$\lambda$. 
\endproclaim

\demo{Proof} First of all, by Theorem 3.4, we can replace ${\Cal U}$ with 
$D(U_\hbar(\b_+))/(\h=\h^*)$. 

Now, in order to construct the functor 
$F$, it is enough to construct a similar functor between 
the the corresponding categories for the algebras 
$D(\b_+)$ and $D(U_\hbar(\b_+))$, i.e. between certain categories of 
dimodules over $\b_-$ and $U_\hbar(\b_-)$. 

But such a functor was constructed in Theorem 4.1. 
Indeed, since all our constructions are compatible 
with the weight decompositions, the functor of Theorem 4.1
restricts to an equivalence $F$ between the categories 
$\Cal O[[\hbar]]$ and $\Cal O_\hbar$. 

The second statement is obvious from the construction. Namely,  
it is easy to see that any formal deformation of a
highest weight $\g$-module to a $U_\hbar(\g)$-module that has the
same character is necessarily a highest weight module.
This fact and the compatibility of $F$ with the weight decomposition 
imply that under $F$, a highest weight module goes to a highest weight
module, a Verma module to a Verma module, and an irreducible
module to an irreducible module (with the same highest weight).
Indeed, the first statement follows from the fact 
that any highest weight module with character equal to the character of
$M(\lambda)$ is isomorphic to $M(\lambda)$, and the second one from
the fact that a formal deformation of an irreducible module is
irreducible. The theorem is proved. 
$\square$ \enddemo

\proclaim{Corollary 4.3} The characters of irreducible highest
weight modules over ${\Cal U}$ are the same as those for
$U(\g)$. 
\endproclaim

\vskip .05in

Corollary 4.3 answers positively question 8.1 from \cite{Dr1}. 

{\bf Remark.} In fact, it is easy to see that Theorem 4.2 implies
a positive answer to Question 8.2 of \cite{Dr1}.
Namely, following \cite{Dr1}, define $I_\beta$ to be the left
ideal in $U(\g)$ generated by elements of weight $\le \beta$. 
We can define a similar ideal $I_\beta^\hbar$ in $U_\hbar(\g)$. 
Then the modules $U(\g)/I_\beta$, $U_\hbar(\g)/I_\beta^\hbar$ 
are in the categories ${\Cal O}$, ${\Cal O}_\hbar$, and we have 
$F(U(\g)/I_\beta)=U_\hbar(\g)/I_\beta^\hbar$, by a deformation
argument. (Indeed, the module $U(\g)/I_\beta$ is generated by a vector $v$
with the only relation $I_\beta v=0$, and $F(U(\g)/I_\beta)$ is a deformation
of $U(\g)/I_\beta$ which has a vector killed by $I_\beta^\hbar$, because
weights are preserved by our construction, so it's
$U_\hbar(\g)/I_\beta^\hbar$). Thus, we have natural isomorphisms 
$\psi_\beta: \text{End}_{U(\g)}(U(\g)/I_\beta)[[\hbar]]\to 
\text{End}_{U_\hbar(\g)}(U_\hbar(\g)/I_\beta^\hbar)$. 
Taking the inverse limit of $\psi_\beta$ with respect to $\beta$
(as the multiplicities of simple roots in $\beta$ go to $+\infty$), we get 
an isomorphism of algebras $\psi: \hat U(\g)[[\hbar]]\to \hat U_\hbar(\g)$, 
required in Question 8.2 of \cite{Dr1}. 

\proclaim{Corollary 4.4} (The Drinfeld-Kohno theorem for
$\g$). Let $k=\Bbb C$. Let $V\in {\Cal O}[[\hbar]]$, and 
$V_q=F(V)$ be its image in $\Cal O_\hbar$. Consider the system of
the Knizhnik-Zamolodchikov 
differential equations with respect to a function 
${\Cal F}(z_1,...,z_n)$ of complex variables $z_1,...,z_n$ 
with values in $ V^{\o n}[\lambda][[\hbar]]$ (the weight subspace of
weight $\lambda$): 
$$
\frac{\partial {\Cal F}}{\partial z_i}=\frac{\hbar}{2\pi {\text{i}}}\sum_{j\ne
i}\frac{\Omega_{ij}{\Cal F}}{z_i-z_j}.
$$
Then the monodromy representation of the braid group $B_n$ for
this equation is isomorphic to the representation 
of $B_n$ on $V_q^{\otimes n}[\lambda]$ defined by the
formula
$$
b_i\to \sigma_i R_{ii+1},
$$
where $b_i$ are generators of the braid group and 
$\sigma_i$ are the permutation of the i-th and (i+1)-th
components. 
\endproclaim

{\bf Remark.} As usual, we identify 
$\pi_1(\Bbb C^n\setminus\{z_i=z_j\}/S_n)$ with $B_n$ by picking 
the reference point $(1,2,...,n)\in \C^n$. 

\demo{Proof} The result follows directly from Theorem 4.2
if we take $\Phi$ to be the Knizhnik-Zamolodchikov associator:
in this case, the two representations are the braid group actions 
on the $n$-th power of two objects in ${\Cal O}$, ${\Cal O}_\hbar$, 
which correspond to each other under the 
braided tensor equivalence $F$.  
$\square$\enddemo

It is easy to generalize these results to the algebra
$\tilde\g'$. Namely, define the categories $\tilde \Cal O[[\hbar]]$
and $\tilde\Cal O$ of representations of 
$\tilde\g'$, $\tilde {\Cal U'}$ similarly to the definition $\Cal
O[[\hbar]]$, $\Cal O_\hbar$. These categories are braided 
in a similar way to $\Cal O[[\hbar]]$, $\Cal O_\hbar$, and we have 

\proclaim{Theorem 4.5} There exists an equivalence of braided 
tensor categories $F: \tilde\Cal O[[\hbar]]\to \tilde\Cal O_\hbar$, 
which is isomorphic to the identity functor at the level of
$\h$-graded $k[[\hbar]]$-modules. This equivalence maps 
the Verma (resp. irreducible) module with highest weight
$\lambda$ to 
the Verma (resp. irreducible) module with highest weight
$\lambda$. 
\endproclaim

\demo{Proof} The proof is the same as the proof of Theorem 4.2, 
using Theorem 3.8.
$\square$\enddemo

{\bf Remark.} Note that in this theorem, we could not 
use $\tilde\g$ instead of $\tilde\g'$, since $\tilde\g$, in general, does not
admit a nondegenerate invariant form, and thus one cannot define 
the element $\Omega$ which is necessary to define the tensor structure. 

\proclaim{Corollary 4.6} The obvious analog of 
Corollary 4.4 is valid if $\g$ is
replaced with $\tilde\g'$. 
\endproclaim

{\bf Remark.} 
We note that Corollary 4.4 for irreducible integrable modules and 
Corollary 4.6 for Verma and contragredient Verma
modules were proved in \cite{V}. 

{\bf 4.3.} Now we want to formulate analytic versions of 
the results of the previous subsection
in which $k=\C$ and $\hbar$ is no longer a formal 
parameter but a complex number. We give such versions in this
subsection. We note that for reader's convenience we do not 
state our results in the maximal possible generality.   

Let us assume for simplicity that the algebra $\g$ 
(in particular, the matrix $A$) is defined over $\Q$
(this is definitely the case for generalized Cartan matrices). 
Let ${\Cal O}_{\Q}$ be the full subcategory of the 
category $\Cal O$ for $\g$
consisting of modules whose weights are defined over $\Q$. 

Let $\hbar\in \C$, $q=e^{\hbar/2}$. By $q^X$ we will always mean
$e^{\hbar X/2}$ . Let $\Cal U_\hbar$ 
be the Drinfeld-Jimbo quantum group, generated 
by $E_i,F_i,q^h$, $h\in \h$, with the usual 
relations, and the relations defined by the kernel 
of the bilinear form $B$. 
Let $\Cal O_{\Q,\hbar}$ be the full subcategory
 of the category $\Cal O$ for $\Cal U_\hbar$, 
consisting of modules whose weights are defined over $\Q$.  

For any $\hbar\in \C$ which is not a nonzero rational multiple of
$\pi {\text{i}}$
(i.e. is such that $q$ is not a nontrivial root of unity), one can define
the structure of a tensor category on both 
${\Cal O}_\Q$ and $ \Cal O_{\Q,\hbar}$. 

Indeed, by standard facts 
about linear ordinary differential equations, the series in
$\hbar$ obtained by restricting the Knizhnik-Zamolodchikov
associator $\Phi$ to a weight subspace in the tensor product
of three objects ${\Cal O}_\Q$ is convergent for small $\hbar$,
and the resulting analytic function continues (in a single-valued
fashion) to the values of $\hbar$ not belonging to $\pi {\text{i}}\Q$
(as for such $\hbar$ the eigenvalues of the operator $\Omega_{ij}$
never differ by a nonzero integer, i.e., no resonances occur). 
This allows us to define a
tensor structure on ${\Cal O}_{\Bbb Q}$ (see also \cite{KL}). 

The structure of a tensor category on ${\Cal O}_{\Q,\hbar}$ comes 
from the Hopf algebra structure on ${\Cal U_\hbar}$. 
Moreover, the first category is braided, with braiding 
$e^{\hbar\Omega/2}$, and the second category is braided 
with braiding defined by the R-matrix, which is well defined 
for generic $\hbar$, i.e., outside of a countable set
(indeed, the $R$-matrix is inverse to the Drinfeld pairing, and 
this pairing is nondegenerate for formal $\hbar$, so 
has countably many zeros for numerical $\hbar$). 

\proclaim{Theorem 4.7} If $\hbar$ is generic
(i.e. outside of a countable set), then   
there exists a 
braided tensor functor $F_\hbar:\Cal O_\Q\to \Cal O_{\Q,\hbar}$, 
which is the identity functor at the level of
$\h$-graded vector spaces, and maps Verma modules to Verma
modules and irreducible modules 
to irreducible modules. 
\endproclaim

\demo{Proof} The theorem is proved similarly to Theorem 
4.2. Namely, consider the functor $F$ constructed in 
Theorem 4.2. One can check directly that for any $V\in \Cal
O_\Q$, the structure maps for $F(V)$ are defined by finite
expressions of the associator $\Phi$ and the braiding
$e^{\hbar\Omega/2}$, 
which implies that they make sense for complex $\hbar\notin \pi
{\text{i}}\Q$.   
$\square$ \enddemo

{\bf Remark.} In the case when $\g$ is a finite dimensional
semisimple Lie algebra, the restriction of the functor 
$F_\hbar$ to the category of finite dimensional modules
is the functor constructed in \cite{KL}. 
The general construction is, essentially, by analogy with \cite{KL}.

\proclaim{Corollary 4.8} For generic $\hbar$, the character of 
the irreducible module $L_\hbar(\lambda)$ over $\Cal U$ with 
highest weight $\lambda\in\h(\Q)$
is the same as the character of the corresponding irreducible module 
$L(\lambda)$ over $\g$. 
\endproclaim

\proclaim{Corollary 4.9} For $V\in \Cal O_\Q$, 
the claim of Corollary 4.4. 
remains valid for generic complex $\hbar$.  
\endproclaim

\demo{Proof} This follows from Theorem 4.7 in the same way as 
Corollary 4.4 follows from Corollary 4.2.   
$\square$\enddemo
 
{\bf Remark.} It is easy to generalize these results 
to the case when $\g$ is replaced with $\tilde\g'$, 
using Theorem 4.5. 

\vskip .05in

Now we would like to make some sharper statements, 
i.e. statements which hold for $\hbar\notin \pi {\text{i}}\Q$. 
To do this, we will assume for simplicity that 
$\g$ is a Kac-Moody algebra. In this case, it is known that the
universal R-matrix is well defined outside of roots of unity, and
that the nilpotent subalgebras of $\Cal U_\hbar$
have the same size as those for $\g$. This allows to
strengthen the above statements (using the same proofs)
as follows. 

\proclaim{Theorem 4.10} If $\hbar\notin \pi \text{i}\Q$, then   
there exists a 
braided tensor functor $F_\hbar:\Cal O_\Q\to \Cal O_{\Q,\hbar}$, 
which is the identity functor at the level of
$\h$-graded vector spaces, and maps Verma modules to Verma
modules and integrable modules 
to integrable modules. 
\endproclaim

{\bf Remark.} We expect that the functor $F_\hbar$ is an
equivalence
if $\hbar\notin \pi {\text{i}}\Q$. 

\proclaim{Corollary 4.11} 
If $V$ is in $\Cal O_\Q$, then the claim of Corollary 4.4. 
remains valid for $\hbar\notin \pi {\text{i}}\Q$.  
\endproclaim

{\bf Remark.} Corollary 4.9 was proved by Drinfeld in the case 
when $\g$ is finite dimensional (\cite{Dr3}). 
Corollary 4.11 for integrable modules 
was proved by Varchenko in \cite{V}.

\end
\Refs
\ref\by [Dr1] V.Drinfeld\paper On some unsolved problems in quantum
group theory \jour Lecture Notes in Math.\vol 1510\pages 1-8\yr
1992\endref

\ref\by [Dr2] V.Drinfeld\paper Quantum groups
\jour Proc. ICM-86 (Berkeley)\vol 1 \yr 1987\pages 798-820\endref

\ref\by [Dr3] V.Drinfeld\paper Quasi-Hopf algebras\jour 
Leningrad Math.J.\vol 1(6)\yr 1990\pages 1419-1457
\endref 

\ref\by [Dr4] V.Drinfeld\paper
On quasitriangular quasi-Hopf algebras and a certain group
closely connected with $Gal(\bar Q/Q)$\jour Leningrad Math. J.\vol
2\yr 1991\pages 829-860\endref

\ref\by [EG] B. Enriquez and N. Geer\paper 
Compatibility of quantization functors 
of Lie bialgebras with duality and doubling operations, 
\yr 2007 \jour math/0707.2337\endref

\ref\by [EK1] P. Etingof and D.  Kazhdan\paper Quantization of
Lie bialgebras, I,
q-alg 9506005\jour Selecta math. \vol 2\issue 1\yr 1996\pages 1-41\endref

\ref\by [EK2] P. Etingof and D. Kazhdan\paper Quantization of Lie bialgebras, 
II, q-alg 9701038
\yr 1998\jour  
Selecta Math.\vol 4\pages 213-231\endref

\ref\by [EK3] P. Etingof and D. Kazhdan\paper Quantization of Lie bialgebras, 
III, q-alg 9610030\yr 1998\jour  
Selecta Math.\vol 4\pages 233-269\endref

\ref\by [EK4] P. Etingof and D. Kazhdan\paper Quantization of Lie bialgebras, 
IV, math.QA/9801043 \yr 2000\jour Selecta Math. \vol 6\issue
1\pages 79-104\endref

\ref\by [EK5] P. Etingof and D. Kazhdan\paper Quantization of Lie bialgebras, 
V, math.QA/9808121 \yr 2000\jour Selecta Math. \vol 6\issue
1\pages 105-130\endref 

\ref\by [J] M. Jimbo\paper 
A q-analog of $U(gl(N+1))$, Hecke algebra, and the
Yang-Baxter equation\jour Lett. Math. Phys. \vol 11 \yr 1986\pages 247-252\endref

\ref\by [K] Kac, V.\book Infinite dimensional Lie algebras \publ 
Cambridge University Press\publaddr Cambridge
\yr 1990\endref

\ref\by [KL] Kazhdan, D., and Lusztig, G.\paper
Tensor structures arising from affine Lie algebras, III\jour 
JAMS\vol 7\issue 2\yr 1994\endref
 
\ref\by [L] Lusztig, G.\book Introduction to quantum groups\publ
Birkh\"auser\publaddr Boston\yr 1994\endref

\ref\by [V] Varchenko, A.\book Multidimensional hypergeometric
functions and representation theory of quantum groups
\publ Adv. Ser. Math. Phys., v.21, World Scientific\publaddr River Edge, NJ
\yr 1995\endref 
\endRefs
